\documentclass[a4paper,11pt,leqno]{article}

\usepackage{amsmath}
\usepackage{amssymb}
\usepackage{amsthm}
\usepackage{bbm}
\usepackage{hyperref}
\usepackage[all]{xy}

\title{A finitely presented subgroup of the automorphism group of a right-angled Artin group}
\author{Emmanuel Toinet}
\date{}

\begin{document}

\maketitle

\begin{abstract}
Let $G_\Gamma$ be a right-angled Artin group. We use geometric methods to compute a presentation of the subgroup $H$ of $Aut(G_\Gamma)$ consisting of the automorphisms that send each generator to a conjugate of itself. This generalizes a result of McCool on basis-conjugating automorphisms of free groups.
\end{abstract}

\section{Introduction}

\newtheorem{definition}{Definition}[section]
\newtheorem{theorem}[definition]{Theorem}
\newtheorem{proposition}[definition]{Proposition}
\newtheorem{lemma}[definition]{Lemma}
\newtheorem{corollary}[definition]{Corollary}
\theoremstyle{remark}
\newtheorem{remark}[definition]{Remark}

\hspace{5mm}A \textit{right-angled Artin group} is a finitely-generated group subject to the relations that some of the generators commute. It can
be described by a simplicial graph $\Gamma$ with one vertex for each generator and one edge for each pair of commuting generators. At the two extremes are the free group $F_n$ of rank $n$ ($\Gamma$ is discrete), and the free abelian group $\mathbbm{Z}^n$ ($\Gamma$ is complete). Right-angled Artin groups are sometimes called \textit{graph groups}, or \textit{free partially commutative groups}. In recent years, right-angled Artin groups have received considerable attention due to the fact that they contain many interesting subgroups, and also because of their actions on CAT(0) cube complexes. The automorphism group of an arbitrary right-angled Artin group is less known (see \cite{D}, \cite{CV}, \cite{Mi}, for example). For a general survey on right-angled Artin groups see \cite{C}.

Automorphisms of right-angled Artin groups were first studied by Servatius in \cite{S}. Drawing on Nielsen automorphisms for free groups, Servatius defined four classes of automorphisms -- consisting of inversions, partial conjugations, transvections, and symmetries (see Section 2) --, and conjectured that they generate $Aut(G_\Gamma)$. Servatius proved his conjecture for some classes of right-angled Artin groups -- for example, when $\Gamma$ is a tree. Thereafter Laurence proved the conjecture for arbitrary right-angled Artin groups in \cite{L}. More recently, Day gave a finite presentation for the automorphism group of a general right-angled Artin group (see \cite{D}), which generalizes the presentation that McCool gave for the automorphism group of a free group (see \cite{M1}).
\\

Let $\Gamma$ = $(\mathcal{V},\mathcal{E})$ be a finite simplicial graph, and let $G_\Gamma$ be the right-angled Artin group associated to $\Gamma$. Recall that $G_\Gamma$ has the presentation:
\begin{center}
$G_\Gamma$ = $\langle$ $\mathcal{V}$ $\mid$ $vw$ = $wv$, $\forall$ \{$v$,$w$\} $\in$ $\mathcal{E}$ $\rangle$.
\end{center} 

\begin{definition}
We say that an automorphism $\varphi$ of $G_\Gamma$ is vertex-conjugating if $\varphi(v)$ is conjugate to $v$ for all $v$ $\in$ $\mathcal{V}$.
\end{definition}

Vertex-conjugating automorphisms were first introduced by Laurence in \cite{L}, where they are called conjugating. They also appear in the recent work of Duncan and Remeslennikov (see \cite{DR}). One of the steps in Laurence's proof of Servatius' conjecture was to show that the set of vertex-conjugating automorphisms coincides with the subgroup $H$ of $Aut(G_\Gamma)$ generated by the partial conjugations (see Section 2 for the definition of a partial conjugation). Let $S$ denote the set of all partial conjugations of $G_\Gamma$. In Section 3, we define a finite set $R$ of relations satisfied by the elements of $S$. Our main result is the following:

\begin{theorem}\label{1.2}
The group $H$ has the presentation $\langle S\mid R\rangle$.
\end{theorem}

In order to prove Theorem \ref{1.2}, we shall construct a finite, connected 2-complex with fundamental group $H$ = $\langle S \mid R \rangle$. Our proof is similar to that of \cite{M3}. In \cite{M3}, McCool gave a finite presentation for the subgroup $H$ of $Aut(F_n)$ consisting of \textit{basis-conjugating} automorphisms. Note that we cannot hope for a generalization of the presentation given in the theorem of \cite{M3} (see \textit{Remark} \ref{3.2} below). Our proof will use the presentation Day gave for the automorphism group of a right-angled Artin group, that will be described in the next section.
\\

I am grateful to Luis Paris, my Ph.D. thesis advisor, for his comments on earlier versions of this work.

\section{Preliminaries}

\hspace{5mm}Let $\Gamma$ = $(\mathcal{V},\mathcal{E})$ be a finite simplicial graph, and let $G_\Gamma$ be the right-angled Artin group associated to $\Gamma$. Let $v$ be a vertex of $\Gamma$. The \textit{link} of $v$, denoted by $lk(v)$, is the subset of $\mathcal{V}$ consisting of all vertices that are adjacent to $v$. The \textit{star} of $v$, denoted by $st(v)$, is $lk(v) \cup \{v\}$. We set $L$ = $\mathcal{V}\cup\mathcal{V}^{-1}$. Let $x$ $\in$ $L$. The \textit{vertex} of $x$, denoted by $v(x)$, is the unique element of $\mathcal{V}\cap\{x,x^{-1}\}$. We set $lk_L(x)$ = $lk(v(x)) \cup lk(v(x))^{-1}$, and $st_L(x)$ = $st(v(x)) \cup st(v(x))^{-1}$.
\\

Let $w$ be a word in $\mathcal{V}\cup\mathcal{V}^{-1}$. The \textit{support} of $w$, denoted by $supp(w)$, is the subset of $\mathcal{V}$ of all vertices $v$ such that $v$ or $v^{-1}$ is a letter of $w$. A word $w$ in $\mathcal{V}\cup\mathcal{V}^{-1}$ is said to be \textit{reduced} if it contains no subwords of the form $vWv^{-1}$ or $v^{-1}Wv$ with $supp(W)$ $\subset$ $star(v)$. For a word $w$ in $\mathcal{V}\cup\mathcal{V}^{-1}$, we denote by $\left|w\right|$ the length of $w$. The \textit{length} of an element $g$ of $G_\Gamma$ is defined to be the minimal length of any word representing $g$. Note that the length of $g$ is equal to the length of any reduced word representing $g$. We say that an element $g$ of $G_\Gamma$ is \textit{cyclically reduced} if it can not be written $vhv^{-1}$ or $v^{-1}hv$ with $v$ $\in$ $\mathcal{V}$, and $\left|g\right|$ = $\left|h\right|$ + 2. By \cite{S}, Proposition 2, every element of $G_\Gamma$ is conjugate to a unique (up to cyclic permutation) cyclically reduced element. The \textit{length} of a conjugacy class is defined to be the minimal length of any of its representative elements. Observe that the length of a conjugacy class is equal to the length of a cyclically reduced element representing it. For an $n$-tuple of conjugacy classes $W$, we define the \textit{length} of $W$, denoted by $\left|W\right|$, as the sum of the lengths of its elements ($n$ $\geq$ 1).
\\

Let $v$, $w$ be vertices of $\Gamma$. We use the notation $v$ $\geq$ $w$ to mean $lk(w)$ $\subset$ $st(v)$. We use the notation $v$ $\sim$ $w$ to mean $v$ $\geq$ $w$ and $w$ $\geq$ $v$.
\\

The Laurence-Servatius generators for $Aut(G_\Gamma)$ are defined as follows:

\begin{description}
\item[Inversions:] Let $v$ $\in$ $\mathcal{V}$. The automorphism $\iota_v$ that sends $v$ to $v^{-1}$ and fixes all other vertices is called an \textit{inversion}.
\item[Partial conjugations:] Let $x$ $\in$ $L$, and let $Y$ be a non-empty union of connected components of $\Gamma \setminus st(v(x))$. The automorphism $c_{x,Y}$ that sends each vertex $y$ in $Y$ to $x^{-1}yx$ and fixes all vertices not in $Y$ is called a \textit{partial conjugation}.
\item[Transvections:] Let $v$, $w$ $\in$ $\mathcal{V}$ be such that $v$ $\geq$ $w$. The automorphism $\tau_{v,w}$ that sends $w$ to $vw$ and fixes all other vertices is called a \textit{transvection}.
\item[Symmetries:] Let $\varphi$ be an automorphism of the graph $\Gamma$. The automorphism $\phi$ given by $\phi(v)$ = $\varphi(v)$  for all $v$ $\in$ $\mathcal{V}$ is called a \textit{symmetry}.
\end{description}

Our aim is to compute a presentation of the subgroup $H$ of $Aut(G_\Gamma)$ generated by the partial conjugations. Our proof will use the fact that partial conjugations are long-range Whitehead automorphisms.
\\

Following \cite{D}, we call \textit{Whitehead automorphism} every automorphism $\alpha$ of one of the following two types:

\begin{description}
\item[Type 1:] $\alpha$ restricted to $\mathcal{V} \cup \mathcal{V}^{-1}$ is a permutation of  $\mathcal{V} \cup \mathcal{V}^{-1}$.
\item[Type 2:] There is an element $a$ $\in$ $L$, called the \textit{multiplier} of $\alpha$, such that $\alpha(a)$ = $a$, and for each $x$ $\in$ $\mathcal{V}$, the element $\alpha(x)$ lies in $\{x,xa,a^{-1}x,\\a^{-1}xa\}$.
\end{description}

One can show that the set of type 1 Whitehead automorphisms is the subgroup of $Aut(G_\Gamma)$ generated by inversions and symmetries.
\\

Following \cite{D}, we say that a Whitehead automorphism $\alpha$ is \textit{long-range} if $\alpha$ is of type 1 or if $\alpha$ is of type 2 and $\alpha$ fixes the vertices of $lk(v(a))$ (where $a$ is the multiplier of $\alpha$).
\\

We denote by $\mathcal{W}$ the set of Whitehead automorphisms, by $\mathcal{W}_1$ the set of Whitehead automorphisms of type 1, and by $\mathcal{W}_2$ the set of Whitehead automorphisms of type 2. We also denote by $\mathcal{W}_\ell$ the set of long-range Whitehead automorphisms.
\\

We use the following notation for type 2 Whitehead automorphisms. Let $A$ be a subset of $L$, and let $a$ $\in$ $L$, such that $a$ $\in$ $A$ and $a^{-1}$ $\notin$ $A$. Provided that it exists, $(A,a)$ denotes the automorphism given by:

\begin{center}
$(A,a)(a)$ = $a$,
\end{center}
and, for all $x$ $\in$ $\mathcal{V}\setminus\{v(a)\}$,
\begin{center}
$(A,a)(x) = \left\{
    \begin{array}{ll}
        x & \mbox{if $x$ $\notin$ $A$ and $x^{-1}$ $\notin$ $A$} \\
        xa & \mbox{if $x$ $\in$ $A$ and $x^{-1}$ $\notin$ $A$} \\
        a^{-1}x & \mbox{if $x$ $\notin$ $A$ and $x^{-1}$ $\in$ $A$} \\
        a^{-1}xa & \mbox{if $x$ $\in$ $A$ and $x^{-1}$ $\in$ $A$}
    \end{array}
\right.$
\end{center}

\vspace{3mm}

If $A$ is a subset of $L$, we set $A^{-1}$ = \{$a^{-1}$ $\mid$ $a$ $\in$ $A$\}. If $A$ and $B$ are subsets of $L$, and $a$ is an element of $L$, we use the notations $A-B$ for $A \setminus B$, $A+B$ for $A \sqcup B$ (if $A \cap B$ = $\emptyset$), $A-a$ for $A\setminus\{a\}$ and $A+a$ for $A\sqcup\{a\}$ (if $a$ $\notin$ $A$). 
\\

The following remark will be of importance in our proof:

\begin{remark}\label{2.1}
Let $x$ $\in$ $L$, and let $Y$ be a non-empty union of connected components of $\Gamma \setminus st(v(x))$. Set $A$ = $Y \cup Y^{-1} \cup \{x\}$, and $a$ = $x$. Then the Whitehead automorphism $(A,a)$ is nothing but the partial conjugation $c_{x,Y}$. In particular, the Whitehead automorphism $(L-lk_L(a)-a^{-1},a)$ is the inner automorphism $\omega_a$ induced by $a$. Note that there is not a unique way to write a partial conjugation as a type 2 Whitehead automorphism. More specifically, if $B$ $\subset$ $lk(v(a))$, then the Whitehead automorphisms $(A,a)$ and $(A+B+B^{-1},a)$ represent the same element of $S$.
\end{remark}

\vspace{4mm}

In \cite{D}, Day proved that $Aut(G_\Gamma)$ is generated by the Whitehead automorphisms, subject to the relations:

\begin{equation}\label{R1}\tag{R1}
(A,a)^{-1} = (A-a+a^{-1},a^{-1}),
\end{equation}
for $(A,a)$ $\in$ $\mathcal{W}_2$.
\begin{equation}\label{R2}\tag{R2}
(A,a)(B,a) = (A \cup B,a),
\end{equation}
for $(A,a)$, $(B,a)$ $\in$ $\mathcal{W}_2$ with $A \cap B$ = \{$a$\}.
\begin{equation}\label{R3}\tag{R3}
(B,b)(A,a)(B,b)^{-1} = (A,a),
\end{equation}
for $(A,a)$, $(B,b)$ $\in$ $\mathcal{W}_2$ such that $a$ $\notin$ $B$, $a^{-1}$ $\notin$ $B$, $b$ $\notin$ $A$, $b^{-1}$ $\notin$ $A$, and at least one of (a) $A \cap B$ = $\emptyset$ or (b) $b$ $\in$ $lk_L(a)$ holds.
\begin{equation}\label{R4}\tag{R4}
(B,b)(A,a)(B,b)^{-1} = (A,a)(B-b+a,a),
\end{equation}
for $(A,a)$, $(B,b)$ $\in$ $\mathcal{W}_2$ such that $a$ $\notin$ $B$, $a^{-1}$ $\notin$ $B$, $b$ $\notin$ $A$, $b^{-1}$ $\in$ $A$, and at least one of (a) $A \cap B$ = $\emptyset$ or (b) $b$ $\in$ $lk_L(a)$ holds.
\begin{equation}\label{R5}\tag{R5}
(A-a+a^{-1},b)(A,a) = (A-b+b^{-1},a)\sigma_{a,b},
\end{equation}
for $(A,a)$ $\in$ $\mathcal{W}_2$, $b$ $\in$ $L$ such that $b$ $\in$ $A$, $b^{-1}$ $\notin$ $A$, $b$ $\neq$ $a$, and $v(b)$ $\sim$ $v(a)$. Here $\sigma_{a,b}$ denotes the type 1 Whitehead automorphism that sends $a$ to $b^{-1}$ and $b$ to $a$, and fixes the other generators.
\begin{equation}\label{R6}\tag{R6}
\sigma (A,a) \sigma^{-1} = (\sigma(A),\sigma(a)),
\end{equation}
for $(A,a)$ $\in$ $\mathcal{W}_2$, and $\sigma$ $\in$ $\mathcal{W}_1$.
\begin{equation}\label{R7}\tag{R7}
\mbox{The entire multiplication table of $\mathcal{W}_1$}
\end{equation}
-- which forms a finite subgroup of $Aut(G_\Gamma)$.
\begin{equation}\label{R8}\tag{R8}
(A,a) = (L-a^{-1},a)(L-A,a^{-1}),
\end{equation}
for $(A,a)$ $\in$ $\mathcal{W}_2$.
\begin{equation}\label{R9}\tag{R9}
(A,a)(L-b^{-1},b)(A,a)^{-1} = (L-b^{-1},b),
\end{equation}
for $(A,a)$ $\in$ $\mathcal{W}_2$, $b$ $\in$ $L$ such that $b$ $\notin$ $A$, $b^{-1}$ $\notin$ $A$.
\begin{equation}\label{R10}\tag{R10}
(A,a)(L-b^{-1},b)(A,a)^{-1} = (L-a^{-1},a)(L-b^{-1},b),
\end{equation}
for $(A,a)$ $\in$ $\mathcal{W}_2$, $b$ $\in$ $L$ such that $b$ $\in$ $A$, $b^{-1}$ $\notin$ $A$, and $b$ $\neq$ $a$.
\\

Note that the relation \eqref{R8} is a direct consequence of the relations \eqref{R1} and \eqref{R2}.
\\

In order to prove Theorem \ref{1.2}, we need to introduce the following technical definitions. 

Let $\alpha$, $\beta$ $\in$ $\mathcal{W}$, and let $W$ be an $n$-tuple of conjugacy classes ($n$ $\geq$ 1). Following \cite{D}, we say that $\beta\alpha$ is a \textit{peak} with respect to $W$ if:
\begin{center}
$\left|\alpha.W\right|$ $\geq$ $\left|W\right|$, \\
$\left|\alpha.W\right|$ $\geq$ $\left|\beta\alpha.W\right|$,
\end{center}
and at least one of these inequalities is strict. 

Let $\alpha_1$,...,$\alpha_k$ $\in$ $\mathcal{W}$ ($k$ $\geq$ 1). We say that $\alpha_i$ is a \textit{peak} of the factorization $\alpha_k\cdots\alpha_1$ with respect to $W$ if 1 $\leq$ $i$ $<$ $k$ and $\alpha_{i+1}\alpha_i$ is a peak with respect to $\alpha_{i-1}\cdots\alpha_i.W$. We say that the factorization $\alpha_k\cdots\alpha_1$ is \textit{peak-reduced} with respect to $W$ if it has no peaks with respect to $W$. The \textit{height} of a peak $\alpha_i$ is $\left|\alpha_i\cdots\alpha_1.W\right|$.

\section{Proof of the main theorem}

\hspace{5mm}In this section, we prove the following:

\begin{theorem}\label{3.1}
The group $H$ has a presentation with generators $c_{x,Y}$, for $x$ $\in$ $L$ and $Y$ a non-empty union of connected components of $\Gamma \setminus st(v(x))$, and relations:
\begin{center}
$(c_{x,Y})^{-1}$ = $c_{x^{-1},Y}$, \\
$c_{x,Y}c_{x,Z}$ = $c_{x,Y \cup Z}$ if $Y \cap Z$ = $\emptyset$, \\
$c_{x,Y}c_{y,Z}$ = $c_{y,Z}c_{x,Y}$ if $v(x)$ $\notin$ $Z$, $v(y)$ $\notin$ $Y$, $x$ $\neq$ $y$, $y^{-1}$, and at least one of $Y \cap Z$ = $\emptyset$ or $y$ $\in$ $lk_L(x)$ holds, \\
$\omega_yc_{x,Y}\omega_y^{-1}$ = $c_{x,Y}$ if $v(y)$ $\notin$ $Y$, $x$ $\neq$ $y$, $y^{-1}$.
\end{center}
\end{theorem}

\textit{Proof}: Our proof is based on arguments developed by McCool in \cite{M2} and \cite{M3} (similar arguments were used in \cite{D}). Recall that $S$ denote the set of partial conjugations. Let $R$ denote the set of relations given in the statement of Theorem \ref{3.1}. We shall construct a finite, connected 2-complex $K$ with fundamental group $H$ = $\langle S \mid R \rangle$. 
\\

We identify a partial conjugation with any of its representatives in $\mathcal{W}_2$ (see Remark \ref{2.1} above). Note that for every $(A,a)$ $\in$ $\mathcal{W}_2$, we have $(A,a)$ $\in$ $S$ if and only if $(A-a)^{-1}$ = $A-a$.
\\

Set $\mathcal{V}$ = \{$v_1$,...,$v_n$\} ($n$ $\geq$ 1). Let $W$ denote the $n$-tuple ($v_1$,...,$v_n$).
\\

The vertices of $K$ will be the set of $n$-tuples $\alpha.W$, where $\alpha$ ranges over the set $\mathcal{W}_1$ of type 1 Whitehead automorphisms. There will be a directed edge $(\alpha.W,\beta\alpha.W;\beta)$ labelled $\beta$ joining $\alpha.W$ to $\beta\alpha.W$ for all $\alpha$, $\beta$ $\in$ $\mathcal{W}_1$. In addition, there will be a loop $(\alpha.W,\alpha.W;(A,a))$ labelled $(A,a)$ at $\alpha.W$ for all $\alpha$ $\in$ $\mathcal{W}_1$, and $(A,a)$ $\in$ $S$. This specifies the 1-skeleton $K^{(1)}$ of $K$.
\\

We shall define the 2-cells of $K$. These 2-cells will derive from the relations \eqref{R1}-\eqref{R10} of \cite{D}. Let $K_1$ be the 2-complex obtained by attaching 2-cells corresponding to the relations \eqref{R7} to $K^{(1)}$. Note that, if $C$ is the 2-complex obtained from $K_1$ by deleting the loops $(\alpha.W,\alpha.W;(A,a))$ for $\alpha$ $\in$ $\mathcal{W}_1$, and $(A,a)$ $\in$ $S$, then $C$ is just the Cayley complex of $\mathcal{W}_1$, and therefore is simply connected.
\\

We now explore the relations \eqref{R1}-\eqref{R5} and \eqref{R8}-\eqref{R10} of \cite{D} to determine which of these will give rise to relations on the elements of $S$.

The relation \eqref{R1} will give rise to the following:
\begin{equation}\label{1}
(A,a)^{-1} = (A-a+a^{-1},a^{-1}),
\end{equation}
for $(A,a)$ $\in$ $S$.

The relation \eqref{R2} will give rise to:
\begin{equation}\label{2}
(A,a)(B,a) = (A \cup B,a),
\end{equation}
for $(A,a)$, $(B,a)$ $\in$ $S$, with $A \cap B$ = \{$a$\}.

The relation \eqref{R3} will give rise to:
\begin{equation}\label{3}
(A,a)(B,b) = (B,b)(A,a),
\end{equation}
for $(A,a)$, $(B,b)$ $\in$ $S$, such that $a$ $\notin$ $B$, $a^{-1}$ $\notin$ $B$, $b$ $\notin$ $A$, and $b^{-1}$ $\notin$ $A$, and at least one of (a) $A \cap B$ = $\emptyset$ or (b) $b$ $\in$ $lk_L(a)$ holds.

From \eqref{R4}, no relations arise. Indeed, suppose that $(A,a)$, $(B,b)$ are in $S$ with $a^{-1}$ $\notin$ $B$, $b$ $\notin$ $A$, and $b^{-1}$
$\in$ $A$. Then $b^{-1}$ = $a$ (because $(A-a)^{-1}$ = $A-a$). But then $a^{-1}$ = $b$ $\in$ $B$ -- leading to a contradiction with our assumption on $a$.

From \eqref{R5}, no relations arise (by the same argument as above).

From \eqref{R8}, we obtain a relation which is a direct consequence of \eqref{1} and \eqref{2}.

The relation \eqref{R9} will give rise to the following:
\begin{equation}\label{4}
(A,a)(L-lk_L(b)-b^{-1},b)(A,a)^{-1} = (L-lk_L(b)-b^{-1},b),
\end{equation}
for $(A,a)$ $\in$ $S$, and $b$ $\in$ $L$ such that $b$ $\notin$ $A$, and $b^{-1}$ $\notin$ $A$.

From \eqref{R10}, no relations arise (by the same argument as above).
\\

We rewrite the relations \eqref{1}-\eqref{4} in the form:
\begin{center}
$\sigma_k^{\varepsilon_k}\cdots\sigma_1^{\varepsilon_1}$ = 1,
\end{center}
where $\sigma_1$,...,$\sigma_k$ $\in$ $S$, and $\varepsilon_1$,...,$\varepsilon_k$ $\in$ \{$-1$,1\}. Let $K_2$ be the 2-complex obtained from $K_1$ by attaching 2-cells corresponding to the relations \eqref{1}-\eqref{4}. Note that the boundary of each of these 2-cells has the form:
\begin{center}
$(\alpha.W,\alpha.W;\sigma_1)^{\varepsilon_1}(\alpha.W,\alpha.W;\sigma_2)^{\varepsilon_2}\cdots(\alpha.W,\alpha.W;\sigma_k)^{\varepsilon_k}$,
\end{center}
for $\alpha$ $\in$ $\mathcal{W}_1$.

Finally, the relations \eqref{R6} will give rise to the following:
\begin{equation}\label{5}
\alpha(A,a)\alpha^{-1} = (\alpha(A),\alpha(a)),
\end{equation}
for $(A,a)$ $\in$ $S$, and $\alpha$ $\in$ $\mathcal{W}_1$. Then $K$ is obtained from $K_2$ by attaching 2-cells corresponding to the relations \eqref{5}. Observe that the boundary of each of these 2-cells has the form:
\begin{center}
$(\beta.W,\beta.W;(\alpha(A),\alpha(a)))^{-1}(\beta.W,\alpha^{-1}\beta.W;\alpha)^{-1}(\alpha^{-1}\beta.W,\alpha^{-1}\beta.W;(A,a))$ \\
$(\alpha^{-1}\beta.W,\beta.W;\alpha)$,
\end{center}
for $\beta$ $\in$ $\mathcal{W}_1$.
\\

It remains to show that $\pi_1(K,W)$ = $H$ = $\langle S \mid R \rangle$.
\\ 

Let $T$ be a maximal tree in $C$. We compute a presentation of $\pi_1(K,W)$ using $T$. There will be a generator $(V_1,V_2;\alpha)$ for each edge $(V_1,V_2;\alpha)$ of $K$.

Since $C$ is simply connected, we have
\begin{equation}\label{6}
(\alpha.W,\beta\alpha.W;\beta) = 1  \hspace{2mm}  \mbox{(in $\pi_1(K,W)$)},
\end{equation}
for all $\alpha$, $\beta$ $\in$ $\mathcal{W}_1$.
\\

Let $\mathcal{P}$ be the set of combinatorial paths in the 1-skeleton $K^{(1)}$ of $K$. We define a map $\widehat{\varphi}$ : $\mathcal{P}$ $\rightarrow$ $Aut(G_\Gamma)$ as follows. For an edge $e$ = $(V_1,V_2;\alpha)$, we set $\widehat{\varphi}(e)$ = $\alpha$, and for a path $p$ = $e_k^{\varepsilon_k}\cdots e_1^{\varepsilon_1}$, we set $\varphi(p)$ = $\widehat{\varphi}(e_k)^{\varepsilon_k}\cdots\widehat{\varphi}(e_1)^{\varepsilon_1}$. Clearly, if $p_1$ and $p_2$ are loops at $W$ such that $p_1$ $\sim$ $p_2$, then $\widehat{\varphi}(p_1)$ = $\widehat{\varphi}(p_2)$. Hence, $\widehat{\varphi}$ induces a map $\varphi$ : $\pi_1(K,W)$ $\rightarrow$ $Aut(G_\Gamma)$. It is easily seen that $\varphi$ is a homomorphism. Then we see from \eqref{6} that $\varphi$ maps $\pi_1(K,W)$ to $H$. It follows immediately from the construction of $K$ that $\varphi$ : $\pi_1(K,W)$ $\rightarrow$ $H$ is surjective. Thus, it suffices to show that $\varphi$ is injective. Let $p$ be a loop at $W$ such that $\varphi(p)$ = 1. We have to show that $p$ $\sim$ 1. Write $p$ = $e_k^{\varepsilon_k}\cdots e_1^{\varepsilon_1}$, where $k$ $\geq$ 1 and $\varepsilon_i$ $\in$ \{$-1$,1\} for all $i$ $\in$ \{1,...,$k$\}. Using the 2-cells arising from the relations \eqref{1}, and the fact that $\mathcal{W}_1^{-1}$ = $\mathcal{W}_1$, we can restrict our attention to the case where $p$ = $e_k\cdots e_1$. Set $\alpha_i$ = $\varphi(e_i)$ for all $i$ $\in$ \{1,...,$k$\}. Note that $\alpha_i$ $\in$ $S \cup \mathcal{W}_1$ $\subset$ $\mathcal{W}_\ell$ for all $i$ $\in$ \{1,...,$k$\}.
\\

Let $Z$ be a tuple containing each conjugacy class of length 2 of $G_\Gamma$, each appearing once.
\\

We shall prove that $p$ $\sim$ $e_l'\cdots e_1'$, such that, if we set $\alpha_i'$ = $\varphi(e_i')$ for all $i$ $\in$ \{1,...,$l$\}, then $\alpha_i'$ $\in$ $\mathcal{W}_1$ or $\alpha_i'$ $\in$ $\mathcal{W}_2\cap Inn(G_\Gamma)$ for each $i$ $\in$ \{1,...,$l$\}. 

First, we examine the case where $\alpha_k\cdots\alpha_1$ is peak-reduced with respect to $Z$. Consider the following sequence of integers:
\begin{center}
$\left|Z\right|$, $\left|\alpha_1.Z\right|$, $\left|\alpha_2\alpha_1.Z\right|$,..., $\left|\alpha_{k-1}\cdots\alpha_1.Z\right|$, $\left|\alpha_k\cdots\alpha_1.Z\right|$ = $\left|Z\right|$.
\end{center}
By Lemma 5.2 in \cite{D}, $\left|Z\right|$  is a minimal element of \{$\left|\alpha.Z\right|$ $\mid$ $\alpha$ $\in$ $\langle\mathcal{W}_\ell\rangle$\}, and therefore is the minimum of the above sequence. On the other hand, since $\alpha_k\cdots\alpha_1$ is peak-reduced with respect to $Z$, there does not exist $i$ $\in$ \{1,...,$k$-1\} such that we have
\begin{center}
$\left|\alpha_{i-1}\cdots\alpha_1.Z\right|$ $\leq$ $\left|\alpha_i\cdots\alpha_1.Z\right|$, \\
$\left|\alpha_{i+1}\cdots\alpha_1.Z\right|$ $\leq$ $\left|\alpha_i\cdots\alpha_1.Z\right|$,
\end{center}
and at least one of these inequalities is strict. Therefore, the above sequence is a constant sequence, and we have
\begin{center}
$\left|\alpha_i\cdots\alpha_1.Z\right|$ = $\left|Z\right|$,
\end{center}
for all $i$ $\in$ \{1,...,$k$\}. We argue by induction on $i$ $\in$ \{1,...,$k$\} to prove that $\alpha_i\cdots\alpha_1.Z$ is a tuple containing each conjugacy class of length 2 of $G_\Gamma$, each appearing once. The result holds for $i$ = 0 by assumption. Suppose that $i$ $\geq$ 1, and that the result holds for $i-1$. Observe that a type 1 Whitehead automorphism does not change the length of a conjugacy class. Thus, we can assume that $\alpha_i$ is a type 2 Whitehead automorphism. Since $\left|\alpha_i\alpha_{i-1}\cdots\alpha_1.Z\right|$ = $\left|\alpha_{i-1}\cdots\alpha_1.Z\right|$, $\alpha_i$ is trivial, or an inner automorphism by \cite{D}, Lemma 5.2. Thus, the result holds for $i$. In this case, $p$ has already the desired form.

We define:
\begin{center}
$h_p$ = max\{$\left|\alpha_i\cdots\alpha_1.Z\right|$ $\mid$ $i$ $\in$ \{0,...,$k$\}\},
\end{center}
and:
\begin{center}
$N_p$ = $|$\{$i$ $\mid$ $i$ $\in$ \{0,...,$k$\} and $\left|\alpha_i\cdots\alpha_1.Z\right|$ = $h_p$\}$|$.
\end{center}
We argue by induction on $h_p$. The base of induction is $\left|Z\right|$ -- the smallest possible value for $h_p$ by \cite{D}, Lemma 5.2. If $h_p$ = $\left|Z\right|$, then the factorization $\alpha_k\cdots\alpha_1$ is peak-reduced and we are done. Thus, we can assume that $h_p$ $>$ $\left|Z\right|$ and that the result has been proved for all loops $p'$ with $h_{p'}$ $<$ $h_p$. Let $i$ $\in$ \{1,...,$k$\} be such that $\alpha_i$ is a peak of height $h_p$. An examination of the proof of Lemma 3.18 in \cite{D} shows that $e_{i+1}e_i$ $\sim$ $f_j\cdots f_1$ such that, if we set $\beta_\kappa$ = $\varphi(f_\kappa)$ for all $\kappa$ $\in$ \{1,...,$j$\}, then:
\begin{equation}\label{7}
\left|\beta_\kappa\cdots\beta_1\alpha_{i-1}\cdots\alpha_1.Z\right| {} < {} \left|\alpha_i\alpha_{i-1}\cdots\alpha_1.Z\right|,
\end{equation}
for all $\kappa$ $\in$ \{1,...,$j-1$\}. Therefore, we get $p$ $\sim$ $e_k\cdots e_{i+2}f_j\cdots f_1e_{i-1}\cdots e_1$ = $p'$, and a new factorization $\alpha_k\cdots\alpha_{i+2}\beta_j\cdots\beta_1\alpha_{i-1}\cdots\alpha_1$. We argue by induction on $N_p$. If $N_p$ = 1, then \eqref{7} implies that $h_{p'}$ $<$ $h_p$ and we can apply the induction hypothesis on $h_p$. If $N_p$ $\geq$ 2, then \eqref{7} implies that $h_{p'}$ = $h_p$ and $N_{p'}$ $<$ $N_p$, and we can apply the induction hypothesis on $N_p$. This completes the induction.
\\

Now, using the 2-cells arising from the relations \eqref{5}, we obtain $p$ $\sim$ $h_s\cdots h_1g_r\cdots g_1$, where, if we set $\gamma_i$ = $\varphi(g_i)$ for all $i$ $\in$ \{1,...,$r$\} and $\delta_j$ = $\varphi(h_j)$ for all $j$ $\in$ \{1,...,$s$\}, then $\delta_i$ $\in$ $\mathcal{W}_1$ for all $i$ $\in$ \{1,...,$s$\} and $\gamma_j$ $\in$ $\mathcal{W}_2\cap Inn(G_\Gamma)$ for all $j$ $\in$ \{1,...,$r$\}. Using \eqref{6}, we obtain $p$ $\sim$ $g_r\cdots g_1$. Set $\mathcal{Z}$ = $\cap_{v \in \mathcal{V}}st(v)$. It follows from Servatius' Centralizer Theorem (see \cite{S}) that the center $Z(G_\Gamma)$ of $G_\Gamma$ is the special subgroup of $G_\Gamma$ generated by $\mathcal{Z}$. Let $\Gamma'$ be the full subgraph of $\Gamma$ spanned by $\mathcal{V}\setminus\mathcal{Z}$. We have $G_{\Gamma'}$ $\simeq$ $Inn(G_\Gamma)$, where the isomorphism is given by $v$ $\mapsto$ $\omega_v$ (see \cite{D}, Lemma 5.3, for example). Write $\gamma_i$ = $(L-lk_L(c_i)-c_i^{-1},c_i)$, where $c_i$ $\in$ $\mathcal{V}\setminus\mathcal{Z}\cup(\mathcal{V}\setminus\mathcal{Z})^{-1}$ ($i$ $\in$ \{1,...,$r$\}). Since $\gamma_r\cdots\gamma_1$ = 1 (in $Inn(G_\Gamma)$), we have $c_r\cdots c_1$ = 1 (in $G_{\Gamma'}$). Therefore $c_r\cdots c_1$ is a product of conjugates of defining relators of $G_\Gamma$. Using the 2-cells corresponding to the relations \eqref{1} and \eqref{3}(b), we deduce that $p$ $\sim$ 1. We conclude that $\varphi$ is injective, and thus $H$ = $\pi_1(K,W)$.
\\

Now, using the 2-cells arising from the relations \eqref{5} (with $\alpha$ = $\beta$), we obtain:
\begin{center}
$(\alpha.W,\alpha.W;(\alpha(A),\alpha(a)))$ = $(\alpha.W,W;\alpha^{-1})(W,W;(A,a))(W,\alpha.W;\alpha)$,
\end{center}
and then, using \eqref{6},
\begin{center}
$(\alpha.W,\alpha.W;(\alpha(A),\alpha(a)))$ = $(W,W;(A,a))$,
\end{center}
for all $\alpha$ $\in$ $\mathcal{W}_1$, and $(A,a)$ $\in$ $S$. It then follows that $H$ is generated by the $(W,W;(A,a))$, for $(A,a)$ $\in$ $S$. We identify $(W,W;(A,a))$ with $(A,a)$ for all $(A,a)$ $\in$ $S$. Any relation in $H$ = $\pi_1(K,W)$ will come from the 2-cells of $K$. Then we see from \eqref{5} that these relations will result from the relations \eqref{1}-\eqref{4} above. It is easily seen that the relations \eqref{1}-\eqref{4} above are equivalent to those of $R$. We have shown that $H$ has the presentation $\langle S \mid R \rangle$.\hfill$\square$
\\

\begin{remark}\label{3.2}
We cannot hope for a generalization of the presentation given in the theorem of \cite{M3}, since, in a general right-angled Artin group, the existence of one-term partial conjugations depends on the existence of domination relations between the vertices of $\Gamma$. (A \textit{one-term partial conjugation} is a partial conjugation of the form $c_{x,\{y\}}$ with $x$ $\geq$ $y$.)
\end{remark}

\begin{flushleft}
Emmanuel Toinet, Institut de Math\'ematiques de Bourgogne, UMR 5584 du CNRS,
Universit\'e de Bourgogne, Facult\'e des Sciences Mirande, 9 avenue Alain Savary, BP 47870, 21078 Dijon Cedex, France \\
E-mail: Emmanuel.Toinet@u-bourgogne.fr
\end{flushleft}


\begin{thebibliography}{9}
\bibitem[Mi]{Mi}
A. Minasyan. Hereditary conjugacy separability of right angled Artin groups and its applications. To appear in Groups, Geometry and Dynamics.
\end{thebibliography}
\end{document}